\documentclass{article}
\usepackage{amssymb,amsmath,enumerate}
\usepackage{a4,latexsym, parskip}
\usepackage{hyperref} 
\hypersetup{pdfauthor=Matthias Sch"utt}
\hypersetup{pdftitle=New examples of modular rigid Calabi-Yau threefolds}

\newcommand{\Phoch}[1]{\mathbb{P}^{#1}}
\newcommand{\cohom}[4][]{\mathrm{H}_{#1}^{#2}(#3,#4)}
\newcommand{\C} [1][]{\mathbb{C}^{#1}}
\newcommand{\Q} [1] []{\mathbb{Q}_{#1}}
\newcommand{\N} [1][] {\mathbb{N}_{#1}}

\newcommand{\F}[1] {\mathbb{F}_{#1}}
\newcommand{\Z}[1][]{\mathbb{Z}^{#1}}
\begin{document}
\setlength{\unitlength}{1cm}

\begin{center}
\LARGE{\textbf{New examples of modular rigid Calabi-Yau threefolds}}
\end{center}

\vspace{0.7cm}

\begin{center}
\large{Matthias Sch\"utt}
\end{center}

\vspace{0.7cm}

\begin{center}
ABSTRACT
\end{center}

The aim of this article is to present five new examples of modular rigid
Calabi-Yau threefolds by giving explicit correspondences to newforms of weight 4 and levels 10,
17, 21, and 73. 

\vspace{0.35cm}  
{\footnotesize

}

\vspace{0.35cm}  

\section{Introduction}

This section introduces the basic concepts of modularity needed to reformulate the modularity
 conjecture from \cite{SY} as a concrete case of the \emph{Fontaine-Mazur conjecture} (cf.
  \cite{FM}).

Recall that a smooth projective threefold $X$ is called \emph{Calabi-Yau} if it satisfies the %%@
following two conditions:

\begin{enumerate}[(i)]

\item the canonical bundle is trivial: $\mathcal{K}_X=\mathcal{O}_X$,
\item $\cohom{1}{X}{\mathcal{O}_X}=\cohom{2}{X}{\mathcal{O}_X}=0.$

\end{enumerate}

We define $X$ to be \emph{rigid} if it has no infinitesimal deformations, i.e.
 $\cohom{1}{X}{\mathcal{T}_X}=0.$
For a threefold with trivial canonical bundle this is equivalent to the vanishing of the
 cohomology group $\cohom{2}{X}{\Omega_X}$ by Serre duality. Hence, the Hodge diamond of a rigid
  Calabi-Yau threefold looks as follows

\[
\begin{matrix}
&&& 1 &&&\\
&& 0 && 0&&\\
& 0 && h^{1,1} && 0 &\\
1 && 0 && 0 && 1\\
& 0 && h^{2,2} && 0 &\\
&& 0 && 0 &&\\
&&& 1 &&&
\end{matrix}
\]

with $h^{1,1}(X)=h^{2,2}(X)=\frac{1}{2}~e(X).$

The main advantage of restricting to rigid Calabi-Yau threefolds lies in the relatively simple
 middle cohomology group $\cohom{3}{X}{\C}$. To introduce the corresponding $\ell$-adic cohomology
  $\cohom[\acute{e}t]{3}{X}{\mathbb{Q}_{\ell}}$ we will further assume $X$ to be defined over $\Q$
   in the following.
Then the natural action of the absolute Galois group Gal$(\bar{\Q}/\Q)$ on
 $\cohom[\acute{e}t]{3}{X}{\mathbb{Q}_{\ell}}$ induces an $\ell$-adic Galois representation

\[
\varrho_X:~~~\text{Gal}(\bar{\Q}/\Q) \rightarrow 
\text{Aut}(\cohom[\acute{e}t]{3}{X}{\mathbb{Q}_{\ell}}) \cong \text{GL}_2(\Q[\ell])
\]
that is well known to be unramified outside the primes of bad reduction of $X$.

We define the $L$-series of $X$ as that one of $\cohom[\acute{e}t]{3}{X}{\mathbb{Q}_{\ell}},$
i.e. as an Euler product

\[
L(X,s)=L(\cohom[\acute{e}t]{3}{X}{\mathbb{Q}_{\ell}},s)=L^{\ast}(s)~\prod_p~
(\text{det}~(\mathbf{1}-\varrho_X(\text{Frob}_p)p^{-s}))^{-1}
\]

where the product runs over the good primes,
while the factor $L^{\ast}(s)$ comes from the bad primes. We shall not need the precise definition
 of $L^{\ast}(s)$ as we are only interested in the good Euler factors of $L(f,s)$.

By linear algebra we have
\[
\text{det}~(\mathbf{1}-\varrho_X(\text{Frob}_p)~T)=1-\text{tr}~ \varrho_X(\text{Frob}_p)~T+
\text{det}~ \varrho_X(\text{Frob}_p)~T^2.
\]
Since furthermore det $\varrho_X$ is known to be the third power of the $\ell$-adic cyclotomic
 character (cf. \cite[p.1]{DM}), the trace $\text{tr} ~\varrho_X(\text{Frob}_p)$ determines the
  $p$th Euler factor. We shall later see that this can be computed via the \emph{Lefschetz fixed
   point formula}. 

We are now able to state the modularity conjecture:

\textbf{The modularity conjecture:}\newline
Any rigid Calabi-Yau threefold $X$ defined over $\Q$ is modular; i.e., its $L$-series $L(X,s)$
coincides up to finitely many Euler factors with the Mellin transform $L(f,s)$ of a modular cusp
form $f$ of weight 4 with respect to $\Gamma_0(N)$, where the level $N$ is only divisible by the %%@
primes of bad reduction.

Quite recently, Dieulefait and Manoharmayum succeeded to prove the modularity conjecture for all
rigid Calabi-Yau threefolds over $\Q$ that have good reduction at 3 and 7 or at 5 and another
suitable prime \cite{DM}. Nevertheless explicit correspondences of rigid Calabi-Yau threefolds
and weight 4 cusp forms are still quite rare and only some few levels are known. The aim of
this article is to construct new examples corresponding to the levels 10, 17, 21, and 73 by
adapting an idea of Chad Schoen\cite{S}.
 
This approach involves the modular elliptic surface of level six, $S_1(6)$, with its natural %%@
projection $pr$
onto $\Phoch{1}:$

\textbf{Theorem.}
{\it
Let $\pi$ be a non-trivial automorphism of\hspace{0.2cm}$\Phoch{1}$, interchanging $0, 1$ and
$\infty$.
Then a small resolution of the twisted self fibre product $(S_1(6),pr)\times_{\Phoch{1}}(S_1(6),
\pi\circ pr)$ is a modular rigid Calabi-Yau threefold,
associated to a newform of weight $4$ and level $10, 17, 21$ or $73$.}

\textbf{Acknowledgement.} I would like to thank K. Hulek and H. Verrill for useful discussions as %%@
well as C. Schoen for helpful comments. 
\section{Construction}

This section recalls Schoen's construction (cf. \cite{S}), starting from two relatively minimal,
regular elliptic surfaces (Y,r), (Y',r') with $r,r'$ surjecting onto the projective line 
$\Phoch{1}$. Let $W$ denote their fibre product $(Y,r) \times_{\Phoch{1}} (Y',r')$. In general,
of course, $W$ will not be smooth; in fact, the singularities are the points $(x,x')$ where $x$
and $x'$ are singular points of the fibres of $(Y,r)$ and $(Y',r')$ over a common cusp $s\in
S''=S\cap S',$ where $S$ and $S'$ denote the images of the singular fibres of $(Y,r)$ and
$(Y',r')$ in  $\Phoch{1}$.

In order to avoid singularities worse than ordinary double points, we are going to assume that all
fibres over $S''$ are either irreducible nodal rational curves or cycles of smooth rational
curves. In Kodaira's notation these are of type $I_b$, where $b>0$ denotes the number of
irreducible components. In this case the canonical sheaf $\omega_W$ of $W$ is locally free. An
easy computation using the adjunction formula and Kodaira's canonical bundle formula shows that
$\omega_W$ is trivial when both elliptic surfaces are rational and admit no multiple fibres (cf. %%@
\cite[p.181]{S}).
Furthermore, by the K\"unneth formula and the Leray spectral sequence the vanishing of
$\cohom{1}{W}{\mathcal{O}_W}$ and $\cohom{2}{W}{\mathcal{O}_W}$ are proven in this situation.

In order to obtain a smooth Calabi-Yau we want to resolve the singularities while preserving the
canonical bundle's triviality. For this purpose we consider a \emph{small resolution} $\hat{W}$
that replaces ordinary double points by $\Phoch{1}$s. However, this is an analytic construction
and the resulting variety is not necessarily projective. Nevertheless, in our examples we will be
able to construct small resolutions as successive blow-ups of closed nonsingular one-codimensional
subvarieties. These global divisors will be components of the singular fibres over $S''$, %%@
isomorphic to $\Phoch{1} \times \Phoch{1}$ (cf. \cite[Lemma(3.1)(i)]{S}).

In order to give a criterion for rigidity Schoen uses different explicit calculations of the Euler
number $e(\hat{W})$ (equaling twice the number of nodes). Thereby he establishes four general
constellations giving rise to rigid Calabi-Yau threefolds \cite[Thm. 7.1]{S}. 

However, many of the possible constructions are ruled out (for our purposes) by the restriction of
being defined over $\Q$.
On the other hand, the varieties we will consider inherit this property from the modular elliptic %%@
surface
$S_1(6)$, since both the twisting automorphism $\pi$ and the surfaces blown up in the small
resolution will be defined over $\Q$.

\section{The modular elliptic surface $S_1(6)$ and its fibre products}

We are interested in the modular elliptic surface $S_1(6)$, i.e. the universal elliptic curve over
$X_1(6)=\overline{\Gamma_1(6)\backslash\mathfrak{S}_1}$ where $\mathfrak{S}_1$ denotes the upper
half plane in $\C$ and $\Gamma_1(6)=\left\{\gamma\in \text{SL}_2(\Z) \mid \gamma \equiv
\begin{pmatrix} 1 & \ast \\ 0 & 1 \end{pmatrix}~\text{mod}~ 6\right\}$ is an arithmetic subgroup
of $\text{SL}_2(\Z)$ operating on $\mathfrak{S}_1$ as Moebius transformations.

As is well known (cf. \cite{B}), a birational model of $S_1(6)$ can be given as the hypersurface
$S\subseteq\Phoch{1}\times\Phoch{2}$ with equation
\[
s(x+y)(y+z)(z+x)+txyz=0.
\]

Resolving the singularities at $((0:1),(1:0:0)), ((0:1),(0:1:0)),$ and\linebreak
$((0:1),(0:0:1))$ we obtain $S_1(6)$ which has the following singular fibres:

\begin{center}
\begin{tabular}{c|c|c}
cusp $(s:t)$ & type of fibre & coordinates of the node(s)  \cr
\hline
$(0:1)$ & $I_6$ & rational\cr
$(1:0)$ & $I_3$ & rational\cr
$(1:1)$ & $I_2$ & irrational; in $\mathbb{Q}(\sqrt{-3})$\cr
$(1:-8)$ & $I_1$ & rational\cr
\end{tabular}
\end{center}

The small resolutions of its self fibre product have been shown to be modular, corresponding to %%@
the unique normalized 
cusp form of weight 4 and level 6,\linebreak
$(\eta(\tau)\eta(2\tau)\eta(3\tau)\eta(6\tau))^2$ (cf. \cite{SY}). In accordance to the modularity %%@
conjecture
this perfectly agrees with the bad reduction at 2 and 3, coming from a change of the types of %%@
singular fibres from $I_3$ to $IV$ and from $I_2$ to $III$, respectively, as the resulting %%@
singularities allow no projective small resolution.

Here, we shall modify this construction. Consider the five non-trivial automorphisms of %%@
$\Phoch{1}$ that interchange the cusps $0, 1$ and $\infty$:

\begin{eqnarray*}
\pi_1: t & \mapsto & 1-t\cr
&&\cr
\pi_2: t & \mapsto & ~~\frac{1}{t}\cr
&&\cr
\pi_3: t & \mapsto & \frac{t}{t-1}\cr
&&\cr
\pi_4: t & \mapsto & \frac{1}{1-t}\cr
&&\cr
\pi_5: t & \mapsto & \frac{t-1}{t}.
\end{eqnarray*}

Then the twisted fibre product $(S_1(6), pr)\times_{\Phoch{1}}(S_1(6), \pi_i\circ pr)$ has only %%@
double points above the three cusps 0, 1, and $\infty$ (since $\pi_i(-8)\neq -8$). According to %%@
Schoen's results it possesses a projective small resolution $\hat{W}^i$ that is indeed a rigid %%@
Calabi-Yau threefold defined over $\Q$.

There are two ways for the $\hat{W}^i$ to obtain bad reduction. One possibility is that the %%@
reduced variety might have no projective small resolution, caused by the degeneration of singular %%@
fibres in characteristic 2 and 3 as in the case of the self fibre product. However, it turns out %%@
that even fibres with only one semi-stable factor allow small projective resolutions: In %%@
characteristic 3, a node of an $I_n \times III$ fibre is given in a formal neighbourhood by %%@
$xy-z(2z+wz-w^2)=0$. This is resolved by blowing up the ideal $(x,z)$ which corresponds to a %%@
component of the fibre (as a global divisor).
Similarly, if $p=2$, we locally have an equation $xy-zw(z+w)=0$ for every node of type $I_n \times %%@
IV$. This is resolved by successively blowing up the ideals $(x,z)$ and $(y,w)$. Again, these %%@
correspond to global divisors. As a consequence, $\hat{W}^i$ has bad reduction at 2 (or 3) if and %%@
only if 0 (or 1) is fixed by $\pi_i$ (such that there is a fibre with both factors of additive %%@
type). 

Other bad primes come from the congruence of $-8$ and its image under $\pi_i$ in some $\F{p}$, %%@
leading to an additional double point for $p\neq 2,3$.
Hence, for the different $\hat{W}^i$ the primes of bad reduction are 17; 3 and 7; 2 and 5; 73; 73. 
In the following we shall associate to the $\hat{W}^i$ newforms of weight 4 whose levels are %%@
exactly the products of these primes, as expected by the modularity conjecture.

\section{Galois representations and Livn\'e's theorem}

Let $f=\sum_{n\in\N} a_n q^n$ be a newform of weight 4 and level $N$. Then its Mellin transform
is given by $L(f,s)=\sum_{n\in\N} a_n n^{-s}$ with argument $s\in\C$. By Atkin-Lehner theory it %%@
has an Euler product expansion with 
good Euler factors $(1-a_p~p^{-s}+p^{3-2s})^{-1}$.

We shall make use of a theorem of Deligne\cite{D} which assigns to $f$ a two-dimensional %%@
$\ell$-adic Galois representation
\[
\varrho_f:  \text{Gal}(\bar{\Q}/\Q) \rightarrow \text{GL}_2(\Q[\ell])
\]
that is unramified outside $\ell$ and the prime divisors of $N$. Deligne's construction implies %%@
the two deep properties
\[
~~\text{tr}~\varrho_f(\text{Frob}_p)=a_p~~\text{and} 
\]\[
\text{det}~\varrho_f(\text{Frob}_p)=p^3~~~~
\]
for all other primes $p$. 

Thereby, the good Euler factors of $L(f,s)$ and $L(\varrho_f,s)$ coincide. This offers us the %%@
opportunity to compare the $\ell$-adic Galois representations $\varrho_{\hat{W}^i}$ and %%@
$\varrho_f$. In view of this approach a powerful tool was developed by Ron Livn\'e  when he recast %%@
work of Faltings and Serre which essentially reduces the proof of the isomorphicity of (the %%@
semi-simplifications of) two Galois representations to checking explicit equality for only %%@
finitely many associated Euler factors.

As the proof of the explicit modularity of $\hat{W}^i$ is based on Livn\'e's theorem, we present %%@
it in the following in a slightly simplified version:

\textbf{Theorem}\cite[Thm. 4.3]{L}
{\it
Let $\varrho_1$, $\varrho_2$ be continuous $2$-dimensional $2$-adic Galois representations. Let %%@
$S$ be a finite set of primes, such that both, $\varrho_1$ and $\varrho_2$, are unramified outside %%@
$S$. Denote by $\Q[S]$ the compositum of all quadratic extensions of $\Q$ which are unramified %%@
outside $S$ and let $T$ consist of primes disjoint from $S$, such that %%@
Gal$(\Q[S]/\Q)=\{\text{Frob}_p:~p\in T\}$.\\
Assume further that

\begin{enumerate}[(i)]
\item  $\text{tr}~\varrho_1(\text{Frob}_p)=\text{tr}~\varrho_2(\text{Frob}_p) $ ~for all $p\in T$, 
\item $\text{det}~\varrho_1(\text{Frob}_p)=\text{det}~\varrho_2(\text{Frob}_p)$ ~for all $p\in T$,
\item  $\text{tr}~\varrho_1\equiv\text{tr}~\varrho_2\equiv 0 $~ {\rm mod} $2$,
\item $\text{det}~\varrho_1\equiv\text{det}~\varrho_2$~ {\rm mod} $2$.
\end{enumerate}

Then $\varrho_1$ and $\varrho_2$ have isomorphic semi-simplifications.}

Since as a consequence the good Euler factors of $\varrho_1$ and $\varrho_2$ coincide,
we will apply this result to the 2-adic Galois representations attached to $\hat{W}^i$ and $f$.

Let $S$ consist of $2$ and all primes of bad reduction of $\hat{W}^i$. In our situation this will %%@
contain the prime divisors of the level $N$. In particular, both representations %%@
$\varrho_{\hat{W}^i}$ and $\varrho_f$ will be unramified outside $S$. Furthermore, %%@
$\Q[S]=\Q(\{\sqrt{-1}\}\cup\{\sqrt{s}: s\in S\})$.

We already saw that
$\text{det}~\varrho_{\hat{W}^i}(\text{Frob}_p)=\text{det}~\varrho_f(\text{Frob}_p)=p^3$ 
for all good primes in our situation. Thus condition $(ii)$ holds. Additionally, Chebotarev's %%@
density theorem ensures condition $(iv)$ and simplifies condition $(iii)$ to
\[
(iii)' ~~\text{tr}~\varrho_{\hat{W}^i}(\text{Frob}_p)\equiv
a_p\equiv 0 ~ \text{mod 2 for all primes} ~p\not\in S.~~~~~~~~~~~~~~~~~~~~~~~~~~~~
\]

Subject to this congruence, the proof of the modularity of $\hat{W}^i$ is therefore reduced to %%@
verifying the identity
\[
(i)' ~~\text{tr}~\varrho_{\hat{W}^i}(\text{Frob}_p)=a_p  ~\text{for all} ~p\in %%@
T~~~~~~~~~~~~~~~~~~~~~~~~~~~~~~~~~~~~~~~~~~~~~~~~~~~~~~~
\] 
for an adequate finite set $T$.

In order to determine such a sufficient set, we use the isomorphism
\begin{eqnarray*}
\Xi: \text{Gal}(\Q[S]/\Q) & \rightarrow & ~~~~(\Z/2\Z)^{\# S+1}\cr
g~~~~~ & \mapsto & \left(\frac{1-\chi_s(g)}{2}\right)_{s\in S\cup\{-1\}}
\end{eqnarray*}
where $\chi_s$ denotes the quadratic Galois character cutting out $\Q(\sqrt{s})$. Since\linebreak
$\chi_s(\text{Frob}_p)=(\frac{s}{p})$ where $(\frac{s}{p})$ denotes the Legendre symbol, we only %%@
have to find a set of (preferably small) primes whose Frobenii map onto 
$(\Z/2\Z)^{\# S+1}$ under the isomorphism $\Xi$.

\section{The traces of $\varrho_{\hat{W}^i}$ and $\varrho_f$}

In general, the operation of Frob$_p^*$ on $\cohom[\acute{e}t]{3}{X}{\Q[\ell]}$ can be quite %%@
difficult to compute. Therefore, we relate it to the easier action on the other cohomology groups %%@
by Grothendieck's version of the Lefschetz fixed point formula
\[
\# ~X (\F{p}) = \sum_{j=0}^{2~\text{dim}\hspace{0.1mm}X} (-1)^j \underbrace{\text{tr} %%@
~(\text{Frob}_p^*~;
\cohom[\acute{e}t]{j}{X}{\Q[\ell]})}_{tr_j(p)}
\]
where $X (\F{p})$ is the set of points of $X$ rational over $\F{p}$ and Frob$_p^*$ denotes the %%@
induced action of the geometric Frobenius morphism. It is worth noticing that %%@
$tr_3(p)=\text{tr}~\varrho_{\hat{W}^i}(\text{Frob}_p)$, since %%@
Frob$_p^*=\varrho_{\hat{W}^i}(\text{Frob}_p)$ on $\cohom[\acute{e}t]{3}{\hat{W}^i}{\Q[\ell]}$ up %%@
to conjugation.

In our cases we have obviously $tr_0(p)=1$ and $tr_1(p)=0$. On the other hand, %%@
$\cohom{2}{\hat{W}^i}{\Z}$ is spanned by divisors which are defined over $\Z$. Indeed, Schoen %%@
shows \cite[p. 190]{S} that the Picard group of $\hat{W}^i$ is spanned by the Picard group of the %%@
generic fibre and the components of the singular fibres of $\hat{W}^i$. Since these are all %%@
defined over $\Z$, we have $tr_2(p)=h^{1,1}p=\frac{e(\hat{W}^i)}{2}~p$ and by Poincar\'e duality
\[tr_3(p)=1+p~(1+p)~\frac{e(\hat{W}^i)}{2}+p^3-\# ~\hat{W^i} (\F{p}).
\]

Counting the points over the finite field $\F{p}$ fibrewise, we obtain the expression
\begin{eqnarray*}
\# ~\hat{W}^i (\F{p})&=&\#~ \hat{W}^i_0 (\F{p})+\#~ \hat{W}^i_1 (\F{p})+
\# ~\hat{W}^i_{\infty} (\F{p})\\
&&+\sum_{k=2}^{p-1} (\# ~S_1(6)_k (\F{p}))
(\# ~S_1(6)_{\pi_i^{-1}(k)} (\F{p})).
\end{eqnarray*}
where the subscripts denote the fibres. Note that the first three summands include the exceptional %%@
lines of the small resolution. However, for the contribution of those lines replacing nodes of %%@
$I_2$-fibres, one has to distinguish whether the coefficients of these nodes are rational over %%@
$\F{p}$, that is, whether $-3$ is a quadratic residue modulo $p$ or not.

While the contribution of the first three summands is easily calculated by hand, 
we need the help of a machine to count the points of (an open part of) the other fibres of %%@
$S_1(6)$. In order to prove modularity we have to search for a newform $f$ with Fourier %%@
coefficients $a_p=tr_3(p)$ for all $p\in T$. For this we use a program of W. Stein \cite{St}.
Having thereby fulfilled condition $(i)'$ from Livn\'e's theorem, we only have to check condition %%@
$(iii)'$ to deduce the identity of all good Euler factors of $L(\hat{W}^i,s)$ and $L(f,s)$.

The evenness of $tr_3(p)$ can be derived from the evenness of 
$\# \hat{W}^i(\F{p})$ for $p\neq 2$. Since the contributions of the fibres with the exceptional %%@
lines are always verified to be even in the next section, this reduces to the evenness of the sum %%@
of the other fibres of $S_1(6)$. However, noticing the symmetry $(1:y:z)\in %%@
S_1(6)_k\Leftrightarrow(1:z:y)\in S_1(6)_k$,
we conclude the evenness of all these fibres ($k=2\hdots p-1$) except for that over -8. Hence, the %%@
whole sum can only be odd if -8 is a fixed point of $\pi_i$ which was explicitly excluded as a %%@
case of bad reduction.

The evenness of the coefficients $a_p$ for all $p\not\in S$ can be established by another method %%@
present in Livn\'e's original paper \cite{L}: Assuming the trace of $\varrho_f$ is not even, %%@
consider the kernel of its mod 2-reduction $\bar{\varrho}_f$. Since by assumption it contains an %%@
element of order 3, the Galois extension of $\Q$, cut out by ker $\bar{\varrho}_f$, must have %%@
Galois group $C_3$ or $S_3$ while being unramified outside 2 and the prime divisors of $N$. The %%@
different extensions possible have been classified by the cubic polynomials they are the splitting %%@
fields of (cf. \cite{J}). We shall see that for each such cubic $h$ we can find a prime $p$ which %%@
produces a contradiction as follows: We choose $p$ such that $h$ is irreducible in $\F{p}$, which %%@
implies that the trace tr $\varrho_f(\text{Frob}_p)=a_p$ is odd, since Frob$_p$ has order 3 in %%@
Gal$(\Q(h)/\Q)$. On the other hand, we use W. Stein's program \cite{St} to compute the Fourier %%@
coefficient $a_p$ and see that it is even, giving the desired contradiction.

\section{The modularity of the twisted fibre products}

We shall now make use of Livn\'e's theorem to establish the modularity of the twisted self fibre %%@
products of $S_1(6)$. This approach demands a 
sufficient set of primes $T$ depending on the ramification set $S$ (in the sense and notation of %%@
section 4). The newform coefficients involved are taken from the tables of W. Stein \cite{St} %%@
while a huge list of the appropriate Galois extensions of $\Q$ was provided by J.Jones \cite{J}.

In the following table we first list Euler number, ramification set and contribution of the fibres %%@
containing the exceptional lines. Afterwards we give the images of the Frobenii Frob$_p$ of some %%@
primes $p$ under the isomorphism $\Xi$, thus proving that they form a sufficient set $T$. We %%@
further give the number of points of $\hat{W}^i$ and the corresponding traces $tr_3(p)$. 

Comparing the traces with the coefficients of a weight 4 newform, we can thus, via Livn\'e's %%@
theorem, establish the modularity of the $\hat{W}^i$.

\subsection{The variety $\hat{W}^1$}

$e(\hat{W}^1)=96$, ~$S=\{2,17\}$, ~
$\#~ \hat{W}^1_0 (\F{p})+\#~ \hat{W}^1_1 (\F{p})+
\# ~\hat{W}^1_{\infty} (\F{p})=48~p~(1+p)$ 

{\scriptsize
\begin{center}
\begin{tabular}{c|c|c|c}
$p$ & $\Xi(\text{Frob}_p)$ & $\#\hat{W}^1(\F{p})$ & $tr_3(p)$\cr
\hline
3 & (1,1,1)& 612 & -8\cr
5 & (0,1,1)& 1560 & 6\cr
7 & (1,0,1)&3060 & -28\cr
13 & (0,1,0)&10992 & -58\cr
19 & (1,1,0)&24984 & 116\cr
41 & (0,0,1)&151920 & -342\cr
47 & (1,0,0)&211824 & 288\cr
89 & (0,0,0)&1090224 & -774
\end{tabular}
\end{center}}

The traces coincide with the Fourier coefficients of the unique normalized cusp form of weight 4 %%@
and level 17, that has only rational coefficients,  %%@
$f(\tau)=q-3q^2-8q^3+q^4+6q^5+-\hdots\in\mathcal{S}_4(\Gamma_0(17))$. The evenness of all its %%@
prime coefficients with the exception of $a_2$ and $a_{17}$ is ensured by the non-existence of an %%@
extension of $\Q$ unramified outside $\{2,17\}$ with Galois group $C_3$ or $S_3$.\newpage

\subsection{The variety $\hat{W}^2$}

$e(\hat{W}^2)=80$,~ $S=\{2,3,7\}$,\newline
$\#~ \hat{W}^2_0 (\F{p})+\#~ \hat{W}^2_1 (\F{p})+
\# ~\hat{W}^2_{\infty} (\F{p})=40~(p^2+p)+
\begin{cases}
~~~0, & \text{if} ~p\equiv 1 ~\text{mod} ~3\\ 4~(p+1), & \text{if} ~p\equiv 2 ~\text{mod} ~3
\end{cases}$

{\scriptsize
\begin{center}
\begin{tabular}{c|c|c|c}
$p$ & $\Xi(\text{Frob}_p)$ & $\#\hat{W}^2(\F{p})$ & $tr_3(p)$\cr
\hline
5 &(0,1,1,1)& 1344 & -18\cr
11 &(1,1,0,1)& 6648 & -36\cr
13 &(0,1,0,1)& 9512 & -34\cr
17 &(0,0,1,1)& 17112 & 42\cr
19 &(1,1,1,0)& 22184 & -124\cr
23 &(1,0,0,1)& 34248 & 0\cr
29 &(0,1,1,0)& 59088 & 102\cr
31 &(1,0,1,0)& 69632 & -160\cr
37 &(0,1,0,0)& 106496 & 398\cr
43 &(1,1,1,1)& 155456 & -268\cr
47 &(1,0,0,0)& 193824 & 240\cr
59 &(1,1,0,0)& 347112 & -132\cr
73 &(0,0,0,1)& 605600 & -502\cr
103 &(1,0,1,1)& 1521152 & 56\cr
137 &(0,0,1,0)& 3329952 & -2358\cr
193 &(0,0,0,0)& 8682704 & 4034\cr
\end{tabular}
\end{center}}

The traces equal the coefficients of a newform\linebreak
$f(\tau)=q-3q^2-3q^3+q^4-18q^5+-\hdots\in\mathcal{S}_4(\Gamma_0(21))^{\text{new}}$. 
There are in total 34 Galois extensions with Galois group $S_3$ or $C_3$ and non-ramification %%@
outside $S$. As these are all splitting fields of certain cubics irreducible modulo some prime %%@
$p\in\{5,11,13,19,23,31\}$,
the evenness of these primes' coefficients guarantees that every coefficient $a_p$ of a good prime %%@
is even.
  
\subsection{The variety $\hat{W}^3$}

$e(\hat{W}^3)=66,~ S=\{2,5\}, ~\#~ \hat{W}^3_0 (\F{p})+\#~ \hat{W}^3_1 (\F{p})+
\# ~\hat{W}^3_{\infty} (\F{p})=33~(p^2+p)$

{\scriptsize
\begin{center}
\begin{tabular}{c|c|c|c}
$p$ & $\Xi(\text{Frob}_p)$ & $\#\hat{W}^3(\F{p})$ & $tr_3(p)$\cr
\hline
3 & (1,1,1) & 468 & -8\cr
7 & (1,0,1) & 2196 & -4\cr
11 & (1,1,0) & 5676 & 12\cr
13 & (0,1,1) & 8262 & -58\cr
17 & (0,0,1) & 14946 & 66\cr
29 & (0,1,0) & 53190 & -90\cr
31 & (1,0,0) & 62376 & 152\cr
41 & (0,0,0) & 126186 & -438
\end{tabular}
\end{center}}

The traces coincide with the coefficients of the unique normalized newform of weight 4 and level %%@
10,
$f(\tau)=q+2q^2-8q^3+4q^4+5q^5-16q^6-4q^7+-\hdots$. Because the only $\Q$-extension unramified %%@
outside $\{2,5\}$ with Galois group $S_3$ ($C_3$ does not occur) is the splitting field of the %%@
polynomial $x^3-x^2+2x+2$, its irreducibility in $\F{3}$
allows us to deduce the evenness of all $f$'s good prime coefficients from that of $a_3$.\newpage

\subsection{The varieties $\hat{W}^4$ and $\hat{W}^5$}

For both Calabi-Yau threefolds $\hat{W}$ the computations turn out the same because $\pi_4$ and %%@
$\pi_5$ induce the same structure as they are inverse to each other:

$e(\hat{W})=72, ~S=\{2,73\},~ \#~ \hat{W}_0 (\F{p})+\#~ \hat{W}_1 (\F{p})+
\# ~\hat{W}_{\infty} (\F{p})=36~(p^2+p)$

{\scriptsize
\begin{center}
\begin{tabular}{c|c|c|c}
$p$ & $\Xi(\text{Frob}_p)$ & $\#\hat{W}(\F{p})$ & $tr_3(p)$\cr
\hline
3 & (1,1,0) & 432 & -8 \cr
5 & (0,1,1) & 1200 & 6\cr
7 & (1,0,1) & 2394 & -34\cr
11 & (1,1,1) & 6078 & 6 \cr
17 & (0,0,1) & 15840 & 90\cr
23 & (1,0,0) & 31980 & 60 \cr
37 & (0,1,0) &101556 &-286\cr
41 & (0,0,0) & 130764 & 150
\end{tabular}
\end{center}}

These traces coincide with those of the unique normalized cusp form of weight 4 and level 73 with %%@
rational Fourier coefficients,
$f(\tau)=q+3q^2-8q^3+q^4+6q^5-24q^6
-34q^7-+\hdots-34q^{13}-+\hdots\in\mathcal{S}_4(\Gamma_0(73)),$ whose prime coefficients outside 2 %%@
and 73 are even, since the three cubics in question (for the splitting fields) are each %%@
irreducible modulo 3 or 13 and both, $a_3$ and $a_{13}$, are even.

\textbf{Theorem.}
{\it
The small resolutions $\hat{W}^i ~(i=1,\hdots,5)$ are modular rigid Calabi-Yau threefolds. The %%@
corresponding modular forms are the newforms of levels $17, 21, 10, 73$, and $73$, given in this %%@
section.}

\vspace{1cm}
Matthias Sch\"utt, \newline
Institut f\"ur Mathematik (C),\newline
Universit\"at Hannover, \newline
Welfengarten 1, \newline
30060 Hannover, Germany,\newline
{\tt schuett@math.uni-hannover.de}

\end{document}